\documentclass[12pt, reqno]{amsart}
\usepackage{amsmath, amsthm, amscd, amsfonts, amssymb, graphicx, color}
\usepackage[bookmarksnumbered, colorlinks, plainpages]{hyperref}
\hypersetup{colorlinks=true,linkcolor=red, anchorcolor=green, citecolor=cyan, urlcolor=red, filecolor=magenta, pdftoolbar=true}

\newtheorem{theorem}{Theorem}[section]

\theoremstyle{definition}

\theoremstyle{remark}
\newtheorem{remark}[theorem]{Remark}
\numberwithin{equation}{section}

\begin{document}

\setcounter{page}{1}

\title[Dense Banach subalgebras]{Dense Banach subalgebras of the null sequence algebra which do not satisfy a differential seminorm condition}

\author[L.B. Schweitzer]{Larry B. Schweitzer$^1$}

\address{$^{1}$Department of Statistics and Biostatistics, California State University East Bay, 25800 Carlos Bee Boulevard, Hayward, CA 94542.}

\email{\textcolor[rgb]{0.00,0.00,0.84}{lschweitzer@horizon.csueastbay.edu;
lsch@svpal.org.}}

\subjclass[2010]{Primary 46L87; Secondary 46H10, 46J10, 46B45, 46K99.}

\keywords{$D_1$-subalgebra, spectral invariance, null sequence algebra, differential structure in $C^\star$-algebras.}

\date{Received: May 12, 2016; Accepted: July 20, 2016.}

\begin{abstract}
We construct dense Banach subalgebras $A$ of the 
null sequence algebra $c_0$ 
which are spectral invariant, but do not
satisfy the $D_1$-condition 
$\|ab \|_A \leq K(\|a\|_{\infty} \|b\|_A + \|a \|_A 
\|b \|_{\infty})$, for all $a, b \in A$.
The sequences in $A$ vanish in a 
skewed manner with respect to an unbounded function
$\sigma\colon {\mathbb N} \rightarrow [1, \infty)$.
\end{abstract} \maketitle

\section{Introduction}
We say that $A$ is a {\it dense Banach subalgebra}
of a $C^\star$-algebra $B$ if $A$ is 
a dense subalgebra of $B$, and 
$A$ is a Banach algebra in some norm $\| \cdot \|_A$,
which is stronger than the restriction to $A$
 of the norm on $B$.
We say that $A$ is {\it spectral invariant}
in $B$ if the quasi-invertible elements of $A$
are precisely those elements of $A$ which are
quasi-invertible in $B$. 
Recall that for $a, b \in A$, 
$a \circ b = a + b - ab$, and $b$ is a 
quasi-inverse for $a$ if and only if 
$a \circ b = b \circ a = 0$ (see \cite[Definition 2.1.1]{Palmer}).
Dense subalgebras of a $C^\star$-algebra 
which are spectral invariant,
or which satisfy a differential seminorm condition such 
as $D_1$, are used
to give differential structure to the $C^\star$-algebra,
and they also have applications in noncommutative
differential geometry 
(see \cite{Connes}\cite{Bhatt}).

Let $c_0$ 
be the $C^\star$-algebra of complex-valued 
vanishing sequences, or null sequences, on the
natural numbers ${\mathbb N}=\{0,1,2,\dots\}$,
with pointwise multiplication and involution.  
For $n \in \mathbb N$, 
let $e_n \colon \mathbb N \rightarrow \{0,1\}$ 
be the unit step function at $n$, where
$e_n(k) = 1$ if $k=n$ and $e_n(k)=0$ 
for $k \in \mathbb N \setminus \{n\}$.
Let $c_f$ denote the linear span of $\{ e_n \}_{n=0}^\infty$,
which is a dense ideal in $c_0$.
The following result shows that it is relatively
easy for a dense subalgebra of $c_0$ to
be spectral invariant.

\begin{theorem}\label{thm:motivator} 
Let $A$ be a dense Banach subalgebra of 
$c_0$.  
Assume that the subalgebra of 
finite support functions $c_f$
is a dense subset of $A$.
Then $A$ is spectral invariant in $c_0$.
\end{theorem}

\begin{proof}
Let $A_{qG}$ denote the group of quasi-invertible 
elements of $A$, with group operation 
$\circ$ and identity $0$.
Let $a \in A \setminus A_{qG}$.
Then $A \circ a$ cannot intersect $A_{qG}$, 
or else $a$ would have a quasi-inverse.
So $\bigl( a + A(1-a)  \bigr) \cap A_{qG}=\emptyset$.
Note that $J= A(1-a)$ is an ideal in $A$, 
and is proper since $J \cap  
\bigl( A_{qG} -a  \bigr) = \emptyset$.
Since $A$ is a Banach algebra, $A_{qG}$ 
and its translate $A_{qG}-a$ are open sets in $A$.
Hence $J$ cannot be dense in $A$, and
$c_f \nsubseteq J$, using the hypothesis.
Since $J$
is a linear space, 
there is some $n_0 \in \mathbb N$ for which
$e_{n_0} \notin J$.
Since $J$ is a $c_f$-module not 
containing $e_{n_0}$, 
every element of $J$ must vanish at $n_0$.
This can only happen if
$a(n_0)=1$.  Hence $a$ is 
not quasi-invertible in $c_0$.
\end{proof}

A dense Banach subalgebra $A$ 
is a {\it $D_1$-subalgebra} of $B$ if 
for some constant $K_A>0$, the {\it $D_1$-condition}
\begin{equation}
\| ab \|_A  \leq  K_A
\bigl\{ \|a\|_A \|b\|_B + \|a\|_B \|b\|_A \bigr\}
\label{eq:D1}
\end{equation}
is satisfied
for all $a, b \in A$,
where $\| \cdot \|_A$ is the norm on $A$ and $\| \cdot \|_B$
is the norm on $B$ (see \cite{KissShul}).
Being a $D_1$-subalgebra implies spectral invariance 
(see \cite[Theorem 5 and Lemma 4]{KissShul}),
which raises the question:
Is every dense Banach subalgebra of the null 
sequence algebra,
which satisfies the hypotheses of Theorem \ref{thm:motivator}, also
$D_1$?  The purpose of the present paper 
is to provide a counterexample.
To this end, we think of $c_0$ 
as ${\mathbb C}^2$-valued sequences
vanishing at infinity, $c_0({\mathbb N}, {\mathbb C}^2)$,
where ${\mathbb C}^2$ denotes the 2-dimensional 
commutative $C^\star$-algebra, 
with coordinatewise multiplication and involution.
We identify the two $C^\star$-algebras
using the isomorphism
$\theta \colon c_0 \cong c_0({\mathbb N}, {\mathbb C}^2)$,
\begin{equation}
\theta(f)(n) = \biggl( f(2n), f(2n+1) \biggr),
\label{eq:ThetaIso}
\end{equation}
for $f \in c_0$, $n \in \mathbb N$.
If $A$ is a dense subalgebra,
$D_1$ is satisfied along each $\mathbb C^2$-summand, because it is
finite dimensional.
In Section 2, we construct 
submultiplicative norms 
  $\| \cdot \|_n$ on the $n$th copy of $\mathbb C^2$,
which make the $D_1$-constants $K_n$
become unbounded as $n$ increases.  In Section 3,
we define the Banach algebra norm $\| \cdot \|_A$ as 
the sup of these $\mathbb C^2$-norms to
construct the counterexample.

\section{Some norms on ${\mathbb C}^2$}
Let $r \in {\mathbb R}$ and let $\sigma$ be a constant
in $[1, \infty)$. 
Let $\| {\vec v} \|_{\max} = \max\{ |x|, |y| \}$,
${\vec v}=(x,y) \in \mathbb C^2$,
 denote the  $C^\star$-norm on $\mathbb C^2$.
Define a seminorm on 
${\mathbb C}^2$ 
by
\begin{equation}
\| {\vec v} \|_{r,\sigma} = 
\|
T_{r, \sigma} {\vec v}
\|_{\max}
 =
\max \bigl\{ |x + ry| , 
{\sigma}|y - rx|\bigr\},
\label{eq:TnormDef}
\end{equation}
where
$T_{r,\sigma}$ is the $(2\times2)$-matrix
\begin{equation}
\begin{pmatrix}
1 & r \\ 
-\sigma r & \sigma
\label{eq:Tdef}
\end{pmatrix}. 
\end{equation}
Note that $\Delta = \det(T_{r,\sigma}) = \sigma(1+r^2)>0$,
so $T_{r,\sigma}$ is invertible, and $\| \cdot \|_{r, \sigma}$
is a norm on $\mathbb C^2$.

We want to find the smallest constant $D > 0$ which
satisfies $\| {\vec v} \|_{\max} \leq
D \| {\vec v} \|_{r, \sigma}$ for all ${\vec v} \in \mathbb C^2$.
Then $D$ is the norm of 
\begin{equation*}
T_{r, \sigma}^{-1} = 
\frac{1}{\Delta}
\begin{pmatrix}
\sigma & -r \\ \sigma r & 1
\end{pmatrix}
\end{equation*}
as an operator on $\mathbb C^2$ with $\max$-norm,
\begin{eqnarray}
D &=& \| T^{-1}_{r, \sigma} \|_\text{op} = 
 \max \big\{ 
\bigl\| ( \sigma, -r) \bigr\|_1,
\bigl\| ( \sigma r, 1) \bigr\|_1
  \bigr\} / \Delta 
\nonumber \\
&=& \frac{1}{ \sigma(1+r^2) } \max\bigl\{ \sigma + |r|, \sigma|r| + 1 \bigr\}
 \nonumber \\
& = &
\frac{1}{ 1+r^2 } \max\bigl\{ 1 + |r|/\sigma, |r| + 1/\sigma \bigr\}
\nonumber \\
& \leq &
\frac{1} { 1+r^2 } \max\bigl\{ 1 + |r|, |r| + 1 \bigr\}
\qquad \mbox{since $\sigma \geq 1$}
\nonumber \\
& = & \frac {1 + |r| }{ 1 + r^2},
\label{eq:D}
\end{eqnarray}
since the dual of the $\max$-norm 
 is the $\ell^1$-norm.
As $r$ ranges over all real numbers, the last expression in \eqref{eq:D}
is bounded by $1.21$.

Next we want to find a constant $C>0$ satisfying
$\| {\vec v}_1 * {\vec v}_2 \|_{r, \sigma}
\leq C \| {\vec v}_1 \|_{r, \sigma} 
\| {\vec v}_2 \|_{r, \sigma}$,
for all ${\vec v}_1, {\vec v}_2 \in \mathbb C^2$, where
${\vec v}_1 * {\vec v}_2$ denotes pointwise multiplication. 
This is equivalent to finding the norm of the operator
${\vec u}_1 \otimes {\vec u}_2 \mapsto T_{r, \sigma} (
T_{r, \sigma}^{-1} {\vec u_1}  *
T_{r, \sigma}^{-1} {\vec u_2} )$ from
$\mathbb C^2 \otimes \mathbb C^2$ to $\mathbb C^2$. 
The operator is $\frac{1} { \Delta^2}$ times
\begin{eqnarray}
& &
\begin{pmatrix}
1 & r \\ -\sigma r & \sigma
\end{pmatrix}
\biggl(
\begin{pmatrix}
\sigma & -r \\ \sigma r & 1
\end{pmatrix}
\begin{pmatrix}
u_{11} \\ u_{12}
\end{pmatrix}
*
\begin{pmatrix}
\sigma & -r \\ \sigma r & 1
\end{pmatrix}
\begin{pmatrix}
u_{21} \\ u_{22}
\end{pmatrix}
\biggr)
\nonumber \\
& = &
\begin{pmatrix}
1 & r \\ -\sigma r & \sigma
\end{pmatrix}
\begin{pmatrix}
(\sigma u_{11} - r u_{12})(\sigma u_{21} - r u_{22}) \\
(\sigma r u_{11} + u_{12})(\sigma r u_{21} + u_{22})
\end{pmatrix}
\nonumber \\
& = &
\begin{pmatrix}
1 & r \\ -\sigma r & \sigma
\end{pmatrix}
\begin{pmatrix}
\sigma^2  & -r \sigma & -r \sigma & r^2 \\
\sigma^2 r^2  & \sigma r & \sigma r & 1 
\end{pmatrix}
\begin{pmatrix}
u_{11}u_{21}   \\ u_{11} u_{22} \\ u_{12} u_{21} \\ u_{12} u_{22}
\end{pmatrix}
\nonumber \\
& = &
\begin{pmatrix}
\sigma^2(1 + r^3)  & \sigma (r^2 - r) & \sigma (r^2-r) & r^2+r \\
\sigma^3(r^2-r)  & \sigma^2 (r^2 +r) & \sigma^2 (r^2 + r) &  \sigma (1-r^3)
\end{pmatrix}
\begin{pmatrix}
u_{11}u_{21}   \\ u_{11} u_{22} \\ u_{12} u_{21} \\ u_{12} u_{22}
\end{pmatrix}.
\nonumber
\end{eqnarray}
So the constant $C$ is
\begin{eqnarray}
C & = & \max \bigl\{ \bigl\|  \bigl( 
\sigma^2(1 + r^3), \sigma (r^2 - r), \sigma (r^2-r),  r^2+r 
\bigr)  \bigr\|_1,
\nonumber \\
& &  \qquad \qquad \qquad
\bigl\|  \bigl(
\sigma^3(r^2-r), \sigma^2 (r^2 +r), \sigma^2 (r^2 + r), \sigma (1-r^3)
\bigr) \bigr\|_1
\bigr\} / \Delta^2
\nonumber \\
  & = &   \frac{ \sigma \max \bigl\{
\frac{|1 + r^3|}{ \sigma} + \frac {2 | r^2 -r|}{ \sigma^2} + 
\frac{|r^2 + r| }{ \sigma^3},
 |r^2 - r| + \frac{2 |r^2 + r|}{ \sigma} + \frac{| 1 - r^3|}{ \sigma^2}
\bigr\} }{
(1+r^2)^2} \nonumber \\
 & \leq &  \frac{ \sigma \max \bigl\{
|1 + r^3| + 2 | r^2 -r| + 
|r^2 + r|,
 |r^2 - r| + 2 |r^2 + r|  + | 1 - r^3| 
\bigr\}}{ 
(1+r^2)^2}, 
\nonumber \\
&&
\label{eq:C}
\end{eqnarray}
where the first step used that the dual of the $\max$-norm
is the $\ell^1$-norm, and 
the last step used that $\sigma \geq 1$. 
As $r$ ranges over all real numbers, the last expression in
\eqref{eq:C} is bounded by $2\sigma$.

It follows that  
 $2 \sigma \| {\vec v}_1 * {\vec v}_2 \|_{r, \sigma}
\leq 
\bigl( 2 \sigma \| {\vec v}_1\|_{r, \sigma}   \bigr)
\bigl( 2 \sigma \| {\vec v}_2 \|_{r, \sigma} \bigr)$, 
for ${\vec v}_1, {\vec v}_2 \in \mathbb C^2$.
By the paragraph before,
$\| {\vec v} \|_{\max} \leq
D \| {\vec v} \|_{r, \sigma} < 
2 \| {\vec v} \|_{r, \sigma} = 
\frac{1}{\sigma} \bigl( 2 \sigma 
 \| {\vec v} \|_{r, \sigma}\bigr) $,
for ${\vec v} \in \mathbb C^2$.

\section{The Banach algebras $A_{r, \sigma}$}
In this section, we pass from the case of a 
single copy of ${\mathbb C}^2$ (see Section 2)
to infinitely many
copies of ${\mathbb C}^2$.
Let  $r \in \mathbb R$ and let
$\sigma$ be any unbounded function from $\mathbb N $ to $[1,\infty)$.
Define an infinite matrix $S_{r, \sigma}$
as the direct sum of $(2 \times 2)$-matrices
\begin{equation*}
S_{r, \sigma} = \bigoplus_{n=0}^\infty
2 \sigma(n) T_{r, \sigma(n)},
\end{equation*}
where each $T_{r, \sigma(n)}$ is defined as in \eqref{eq:Tdef}.
Let $A_{r, \sigma} = \{\, f \in c_0 \, |\,
S_{r,\sigma} \theta f \in c_0({\mathbb N}, {\mathbb C}^2) \, \}$,
where $\theta$ was defined in the Introduction (see \eqref{eq:ThetaIso}).
By \cite[Theorem 4.3.1]{Wilansky},
$A_{r, \sigma}$ is complete in the norm
\begin{equation}
\| f\|_{r, \sigma}
\, = \, \|S_{r,\sigma} \theta f \|_{c_0({\mathbb N},{\mathbb C}^2)}
\, = \, \sup_{n=0}^\infty  \, \bigl\{ 2 \sigma(n)  
\bigl\|\bigl(f(2n), f(2n+1)\bigr)\bigr\|_{r, \sigma(n)} \bigr\},
\label{eq:AsigNorm}
\end{equation}
where each $\| \cdot \|_{r, \sigma(n)}$ is defined as in \eqref{eq:TnormDef},
and where by the final remarks of Section 2, $\| f g \|_{r, \sigma} \leq
\| f \|_{r, \sigma} \| g \|_{r, \sigma}$  
and $\| f \|_\infty \leq \| f \|_{r, \sigma}$, 
for $f, g \in A_{r,\sigma}$.
Since  $S_{r, \sigma} \theta f \in c_0({\mathbb N}, {\mathbb C}^2)$
for $f \in A_{r, \sigma}$, then
\begin{equation*}
0 = 
\lim_{n\rightarrow \infty} \, 
\bigl\| \bigl( S_{r, \sigma} \theta f \bigr) (n) \bigr\|_{{\mathbb C}^2}
\,
=
\lim_{n\rightarrow \infty}  \,  2 \sigma(n)  
\bigl\|\bigl(f(2n), f(2n+1)\bigr)\bigr\|_{r, \sigma(n)}. 
\end{equation*}
For $\epsilon>0$, let $N_{\epsilon}$ be large enough
so that  
the argument of the limit
is smaller than $\epsilon$ if $n > N_{\epsilon}$.
Then for $n \geq N_{\epsilon}$,
\begin{eqnarray}
\bigl\| f - \bigl(f(0), \dots, f(2n+1), 0,0, \dots\bigr) \bigr\|_{r,\sigma} &=&
\bigl\| \bigl(\underbrace{0, \dots, 0}_{\text{$2n+2$ zeros}}, 
f(2n+2), \dots\bigr) \bigr\|_{r,\sigma}  
\nonumber \\
 &=& \sup_{k >n} 2 \sigma(k) \bigl\|\bigl( f(2k),
 f(2k+1)\bigr) \bigr\|_{r, \sigma(k)} 
\nonumber \\
 & <& \epsilon,
\nonumber
\end{eqnarray}
using the definition of $\| \cdot \|_{r,\sigma}$
\eqref{eq:AsigNorm}.
It follows that $c_f$ is dense in $A_{r, \sigma}$, 
and we can apply
Theorem 1.1 to see that $A_{r, \sigma}$ is spectral invariant 
in $c_0$.

\begin{theorem}
 The dense Banach subalgebra $A_{r, \sigma}$ 
of $c_0$ is not a $D_1$-subalgebra of $c_0$ 
for $r \in \mathbb R \setminus \{ 0, 1\}$.
\end{theorem}

\begin{proof}
For $n \in \mathbb N$, define $a_n\in A_{r, \sigma}$ by 
\begin{equation*}
a_n = \biggl(\underbrace{0,0, \dots 0,0}_{\text{$2n$ zeros}},
 1,r, 0,0,
	\dots \biggr), 
\end{equation*}
where $a_n(2n)=1$ and $a_n(2n+1)=r$ are the only nonzero components.
Then 
\begin{eqnarray}
\| a_n \|_{r, \sigma} &=& 2\sigma(n)(1+r^2), \nonumber \\
\| a_n^2 \|_{r, \sigma} &=&
2\sigma(n)
\max\bigl( |1+r^3| , \sigma(n)|r^2 - r| \bigr),\nonumber \\
  \| a_n \|_{\infty} 
&=& \max(1, |r|). 
\nonumber
\end{eqnarray}
If $K>0$ were a constant satisfying the  
$D_1$-condition \eqref{eq:D1} for 
$A_{r, \sigma}$ in $c_0$, then 
\begin{eqnarray}
K & \ge  &
\frac{\|a_n^2\|_{r, \sigma}}{2\| a_n\|_{\infty} \| a_n 
\|_{r, \sigma}} = \frac{\max\bigl(|1+r^3| , 
\sigma(n)|r^2 - r|\bigr) }{ 2(1+r^2)\max(1, |r|) }
\nonumber \\
& \geq & \frac{\sigma(n)|r^2-
r| }{ 2(1+r^2)\max(1, |r|) }
\nonumber
\end{eqnarray}
 must hold for each 
$n$ for which $\sigma(n) |r^2 - r| \geq |1 + r^3|$.
No such constant 
$K$ can exist if $r \not= 1$ 
or $0$, 
since $\sigma$ is unbounded.
\end{proof}

\begin{remark}
Note that $A_{r,\sigma}$ is a Banach $\star$-algebra.
The norms defined on $\mathbb C^2$ \eqref{eq:TnormDef}
and the norm 
on $A_{r, \sigma}$
\eqref{eq:AsigNorm} 
are both left unchanged by the $\star$-operation of
pointwise complex-conjugation.
\end{remark}

\begin{remark}
In the cases $r=0$ and $r=1$,
it can be shown that 
$A_{r, \sigma}$ is a $D_1$-subalgebra of
$c_0$.  Further, $A_{r, \sigma}$
 is a dense Banach 
ideal in $c_0$ if $r=0$, but not if $r=1$.
\end{remark}

\bibliographystyle{amsplain}

\end{document}